\documentclass[12pt]{amsart}
\usepackage{amsmath,amsthm,amsfonts,amssymb,latexsym,enumerate,xcolor}
\usepackage[pagebackref]{hyperref}
\usepackage[english]{babel}

\newcommand{\Aut}{{{\operatorname{Aut}}}}

\newcommand{\GL}{\operatorname{GL}}

\newcommand{\PSL}{\operatorname{PSL}}

\newcommand{\PSp}{\operatorname{PSp}}

\newcommand{\Fa}{\operatorname{F^{\ast}}}

\newcommand{\Fi}{\operatorname{F}}

\newcommand{\fpr}{\operatorname{fpr}}

\newcommand{\dd}{\operatorname{d}}
\newcommand{\Pro}{\operatorname{Pr}}

\newcommand{\Sym}{\operatorname{Sym}}

\newtheorem{thm}{Theorem}[section]
\newtheorem{lem}[thm]{Lemma}

\newtheorem{pro}[thm]{Proposition}
\newtheorem{cor}[thm]{Corollary}

\newtheorem*{thmA}{Theorem A}
\newtheorem*{conA'}{Conjecture A'}
\newtheorem*{thmB}{Theorem B}
\newtheorem*{thmC}{Theorem C}

\theoremstyle{definition}

\numberwithin{equation}{section}

\begin{document}





\title[On the commuting probability of $\pi$-elements in finite groups]{On the commuting probability of $\pi$-elements in finite groups}

\author{Juan Mart\'{\i}nez}
\address{Departament de Matem\`atiques, Universitat de Val\`encia, 46100
  Burjassot, Val\`encia, Spain}
\email{Juan.Martinez-Madrid@uv.es}

\thanks{Research supported by Generalitat Valenciana CIAICO/2021/163 and CIACIF/2021/228. The author would like to thank T. C. Burness  and A. Moretó  for many helpful conversations on this paper. Finally, the author thanks the
referee for the careful reading of the paper and helpful comments.}




\begin{abstract}
Let $G$ be a finite group, let $\pi$ be a set of primes and let $p$ be the smallest prime in $\pi$. In this work, we prove that $G$ possesses a normal and abelian Hall $\pi$-subgroup if and only if  the probability that two random $\pi$-elements of $G$ commute is larger than $\frac{p^2+p-1}{p^3}$. We also prove that if $x$ is a $\pi$-element not lying in $O_{\pi}(G)$, then the proportion of $\pi$-elements commuting with $x$ is at most $1/p$.

\end{abstract}

\maketitle



\section{Introduction}

Let $G$ be a finite group. We  define 
$$\Pro(G)=\frac{|\{(x,y)\in G\times G|xy=yx\}|}{|G|^{2}}$$
to be the probability that two random elements of $G$ commute. It was proved in \cite{Gustafson} that this probability coincides with $k(G)/|G|$, where $k(G)$ denotes the number of conjugacy classes of $G$. A result often attributed to Gustafson asserts that  if $\Pro(G)>\frac{5}{8}$, then $G$ is abelian. In fact,  it had been proved before in \cite{J} that  if $p$ is the smallest prime dividing $|G|$ and $\Pro(G)>\frac{p^2+p-1}{p^3}$, then $G$ is abelian (a proof of this result  can also be found in \cite{Orig}). There exists a large literature on the commuting probability, see for example \cite{GR,N}.


 There exist also local versions of  $\Pro(G)$ which provide  information on the structure of Sylow subgroups and Hall subgroups of the group. As a first example, we have  $d_{\pi}(G)$, which was introduced by A. Maróti and H. N. Nguyen \cite{MN}. Given a set of primes $\pi$ and a finite group $G$,  $d_{\pi}(G)$ is defined as 
 $$\dd_{\pi}(G)=\frac{k_{\pi}(G)}{|G|_{\pi}},$$
where $k_{\pi}(G)$ denotes the number of conjugacy classes of $G$ whose elements are $\pi$-elements and $|G|_{\pi}$ denotes the $\pi$-part of $G$. It was proved in \cite{MN} that if $\dd_{\pi}(G)>\frac{5}{8}$, then $G$ possesses an abelian Hall $\pi$-subgroup. This result was improved by the main result of \cite{HMM} which asserts that if $p$ is the smallest prime in $\pi$ and $\dd_{\pi}(G)>\frac{p^{2}+p-1}{p^{3}}$, then $G$  possesses an abelian Hall $\pi$-subgroup.

In this paper, we will  study  a different local invariant. Given a set of primes $\pi$ and a finite group $G$,  we  define $\Pro_{\pi}(G)$ as the probability that two randomly chosen  $\pi$-elements of $G$ commute, that is
$$\Pro_{\pi}(G)=\frac{|\{(x,y)|x,y\in G_{\pi}, xy=yx\}|}{|G_{\pi}|^{2}},$$
where $G_{\pi}$ is the set of all $\pi$-elements in $G$. T. C. Burness, R. M. Guralnick, A. Moretó and G. Navarro \cite{Orig} introduced $\Pro_{p}(G)$, where $p$ is a prime. They proved \cite[Theorem A]{Orig} that $G$ possesses a normal and abelian Sylow $p$-subgroup if and only if $\Pro_{p}(G)>\frac{p^{2}+p-1}{p^{3}}$. In fact, this result was a direct consequence of \cite[Theorem C]{Orig}, which asserts that if $x\in G_p\setminus O_p(G)$, then $x$ commutes with at most $1/p$ of the $p$-elements of $G$.

Our goal is to extend  these results for a general set of primes $\pi$.

\begin{thmA}
Let $G$ be  a finite group, let $\pi$ be a set of primes and let $p$ be the smallest prime in $\pi$. If $x \in G_{\pi}\setminus O_{\pi}(G)$, then 
$$\frac{|C_{G}(x)_{\pi}|}{|G_{\pi}|} \leq \frac{1}{p}.$$
\end{thmA}

Let $\Fa(G)$ be the generalized Fitting subgroup of $G$.  In the case when $\Fa(G)$ is a $\pi'$-group, it is possible to prove a more general version of Theorem A.

\begin{thmB}
Let $\pi$ be a set of primes and let $G$ be  a finite group such that $\Fa(G)$ is a $\pi'$-group. If $x$ is a non-trivial $\pi$-element of $G$ and $q$ is the largest prime dividing $o(x)$, then 
$$\frac{|C_{G}(x)_{\pi}|}{|G_{\pi}|} \leq \frac{1}{q}.$$
\end{thmB}


We remark that Theorem B is independent of Theorem A and its proof does not depend on the classification of finite simple groups. It remains as an open question whether the conclusion of Theorem B holds without assuming that  $\Fa(G)$ is a $\pi'$-group. We do not know if Theorem A holds by replacing the smallest prime in $\pi$ by the smallest prime dividing $o(x)$.

 From Theorem A, we deduce an analogue version of \cite[Theorem A]{Orig} for a general set of primes $\pi$.

\begin{thmC}
Let $G$ be a finite group, let $\pi$ be  a set of primes and let $p$ be the smallest member of $\pi$. The group $G$ has a normal and abelian Hall $\pi$-subgroup if and only  $\Pro_{\pi}(G)> \frac{p^{2}+p-1}{p^{3}}$.
\end{thmC}

We observe  that for every   $p>3$ the group $\PSL(2,p)$ is  simple and  $ \Pr_p(\PSL(2,p))=\frac{p^2+p-1}{p^3}$. Therefore,  \cite[Theorem A]{Orig} and Theorem C are best possible for every $p>3$.


The paper is organized as follows. In Section \ref{Preliminary} we will present some preliminary results. In Sections \ref{TheoremB} and \ref{TheoremA} we will prove Theorems B and A, respectively. Finally, in Section \ref{TheoremC} we will use Theorem A to prove Theorem C.

We remark that the proofs of the main results  of this paper are a refinement of the arguments of \cite{Orig}.

\section{Preliminary results}\label{Preliminary}

In this section we collect the preliminary results that we will use later. Our first result relates $\Pr_{\pi}(G)$ with $\Pr_{\pi}(G/N)$ for $N$,  a normal  $\pi$-subgroup of $G$.

\begin{lem}\label{gruponorm}
Let $\pi$ be a set of primes, let $N$ be a normal $\pi$-subgroup of $G$ and let $x\in G_{\pi}$. Then
\begin{itemize}
\item [(i)] $\frac{|C_{G}(x)_{\pi}|}{|G_{\pi}|}\leq \frac{|C_{G/N}(xN)_{\pi}|}{|(G/N)_{\pi}|}$.

\item [(ii)] $\Pro_{\pi}(G)\leq \Pro_{\pi}(G/N)$.
\end{itemize}
\begin{proof}
Let $\overline{\cdot}$ denote the quotient by $N$.

First, we express $\overline{G}_{\pi}$ in terms of $G_{\pi}$. Since $N$ is a normal $\pi$-subgroup, we have that $\overline{x}$ is a $\pi$-element if and only if $x$ is a $\pi$-element. Thus $\overline{G}_{\pi}=\overline{G_{\pi}}$. Since $\{xN|x \in G_{\pi}\}$ is a partition of $G_{\pi}$, we deduce that $|\overline{G}_{\pi}|=\frac{|G_{\pi}|}{|N|}$.


In order to prove (i), it suffices to show that $|C_G(x)_{\pi}|\leq |N||C_{\overline{G}}(\overline{x})_{\pi}|$. It is easy to see that if $y \in C_G(x)_{\pi}$, then $\overline{y}\in C_{\overline{G}}(\overline{x})_{\pi}$. It follows that
$$|C_{\overline{G}}(\overline{x})_{\pi}|\geq |\{\overline{y}|y \in C_G(x)_{\pi}\}|=\frac{|\{yt|y \in C_G(x)_{\pi},t \in N\}|}{|N|}\geq \frac{|C_G(x)_{\pi}|}{|N|},$$
which gives (i).

Now, let us prove (ii). Applying (i) and that $|G_{\pi}|=|N||\overline{G}_{\pi}|$, we obtain 
$$\Pro_{\pi}(G)=\frac{1}{|G_{\pi}|}\sum_{x \in G_{\pi}}\frac{|C_{G}(x)_{\pi}|}{|G_{\pi}|}\leq \frac{1}{|\overline{G}_{\pi}|}\cdot \frac{1}{|N|}\sum_{x \in G_{\pi}} \frac{|C_{\overline{G}}(\overline{x})_{\pi}|}{|\overline{G}_{\pi}|}.$$ 
Now, we observe that 
$$\sum_{x \in G_{\pi}} \frac{|C_{\overline{G}}(\overline{x})_{\pi}|}{|\overline{G}_{\pi}|}=\sum_{\overline{y} \in \overline{G}_{\pi}}\sum_{n\in N} \frac{|C_{\overline{G}}(\overline{y} )_{\pi}|}{|\overline{G}_{\pi}|}=|N|\sum_{\overline{y}  \in \overline{G}_{\pi}}\frac{|C_{\overline{G}}(\overline{y} )_{\pi}|}{|\overline{G}_{\pi}|}.$$
Thus, we deduce that
$$\Pro_{\pi}(G)\leq \frac{1}{|\overline{G}_{\pi}|}\sum_{\overline{y}  \in \overline{G}_{\pi}}\frac{|C_{\overline{G}}(\overline{y} )_{\pi}|}{|\overline{G}_{\pi}|}=\Pro_{\pi}(\overline{G})$$
and (ii) follows.
\end{proof}
\end{lem}

 Now, we present some results on coprime actions that we will use in the proof of Theorem B.

\begin{lem}[Lemma 2.3 of \cite{Orig}]\label{inequality}
Let $P$ be a $p$-group acting coprimely on a finite group $K$ and let $L$ be a $P$-invariant subgroup of $K$. If 
$$\frac{|C_{K}(P):C_{L}(P)|}{|K:L|}<1,$$
then

$$\frac{|C_{K}(P):C_{L}(P)|}{|K:L|}\leq \frac{1}{p+1}.$$
\end{lem}

\begin{lem}[Lemma 2.4 of \cite{Orig}]\label{conj}
Let $G$ be a finite group, let $x,y \in G$ and let $K\leq G$ such that $x$ and $y$ normalize $K$, $Kx=Ky$ and $(|K|,o(x))=1=(|K|,o(y))$. Then, $x=y^{k}$ for some $k \in K$.
\end{lem}

\begin{cor}\label{res}
Let $G$ be a finite group and let $\pi$ be a set of primes. If $y \in G$ is a $\pi$-element and $K\leq G$  is a $\pi'$-subgroup of $G$ normalized by $y$, then $y^{K}=(Ky)_{\pi}$.
\begin{proof}
Clearly, $y^K$ is contained in $(Ky)_{\pi}$. Now, let $z\in (Ky)_{\pi}$. We know that $Ky=Kz$ and, since $K$ is a $\pi'$-group and $y$ is a $\pi$-element, we have that $(|K|,o(y))=1=(|K|,o(z))$. Thus, by Lemma \ref{conj}, we deduce that $z$ and $y$ must be $K$-conjugate.
\end{proof}
\end{cor}

Finally, we introduce the fixed point ratio, which will be the key for the proof of Theorem A.

Let $G$ be a finite group acting on a finite set $\Omega$. Given $z \in G$, we define the fixed point ratio of $z$ as the proportion of elements in $\Omega$ fixed by $z$. We will denote it by $\fpr(z,\Omega)$.

It is important to remark that given $x$, a $\pi$-element of $G$, the proportion $$\frac{|C_G(x)_{\pi}|}{|G_{\pi}|}$$ is exactly the fixed point ratio of $x$, where we consider the action of $G$ on $G_{\pi}$ by conjugation.

We will need the following results from \cite{FPR} on the fixed point ratios of primitive permutation groups.

\begin{thm}[Theorem 1 of \cite{FPR}]\label{fpr1}
Let $G\leq \Sym(\Omega)$ be a finite primitive permutation group with socle $J$ and
point stabilizer $H$. If $z \in G$ is an element of prime order $p$, then either
$$\fpr(z,\Omega)\leq \frac{1}{p+1},$$
or one of the following holds:

\begin{itemize}

\item [(i)] $G$ is almost simple and one of the following holds:

\begin{itemize}
\item [a)] $G\in \{S_n, A_n \}$ acting on $k$-element subsets of $\{1,\ldots,n\}$ with $1\leq k <n/2$.

\item [b)]  $G$ is classical in a subspace action and    $(G,H,z,\fpr(z,\Omega))$ is   listed in Table 6 of \cite{FPR}.

\item [c)]  $G=S_n$, $H=S_{n/2}\wr C_2$ and $x$ is a transposition.

\item [d)]  $G=M_{22}:2$, $H=L_3(4).2_2$ and $x=2B$.
\end{itemize}

\item [(ii)] $G$ is an affine group, $J = (C_p)^d$,  $z \in \GL_d(p)$ is a transvection and $\\ \fpr(z, \Omega) = \frac{1}{p}$.

\item [(iii)] $G \leq A\wr S_t$ is a product type group with its product action on $\Omega=\Gamma^t$ and $z \in A^t\cap G$, where $A \leq \Sym(\Gamma)$ is one of the almost simple primitive groups in part (i).
\end{itemize}
\end{thm}

\begin{cor}[Corollary 3 of \cite{FPR}]\label{fpr2}
Let $G\leq \Sym(\Omega)$ be a finite almost simple primitive permutation group with socle $J$ and
point stabilizer $H$. If $z \in G$ has prime order $p$. Then either
$$\fpr(z,\Omega)\leq \frac{1}{p},$$
or one of the following holds:

\begin{itemize}

\item [(i)] $J= A_n$ and $\Omega$ is the set of $k$-element subsets of $\{1,\ldots,n\}$ with $1\leq k <n/2$.

\item [(ii)] $G$ is a classical group in a subspace action and $(G,J,p,x,\fpr(x))$ is one of the possibilities recorded in Table 1 of \cite{FPR}.
\end{itemize}
\end{cor}

\section{Proof of Theorem B}\label{TheoremB}

In this section, we prove Theorem B. The following results give information on the proportion of $\pi$-elements which commute with a fixed $\pi$-element when the group possesses a normal and non-trivial $\pi'$-group.

\begin{pro}\label{first}
Let $G$ be a finite group, let $\pi$ be a non-empty set of primes and  let $K$  be a normal $\pi'$-subgroup of $G$. Assume that $x \in G$ is an element of order $p$, where $p \in \pi$, such that $[x,K]=K$. If $y \in G$ is a $\pi$-element which commutes with $x$, then the following holds.
\begin{itemize}
    \item [(i)] If $[y,K]\not=1$, then $$\frac{|\{z \in(Ky)_{\pi}|[x,z]=1 \}|}{|(Ky)_{\pi}|}\leq \frac{1}{p+1}.$$
    
    \item [(ii)] If $L=\langle K,x \rangle$, then $$\frac{|\{z \in(Ly)_{\pi}|[x,z]=1 \}|}{|(Ly)_{\pi}|}\leq\frac{1}{p}.$$
\end{itemize}
\begin{proof}
We begin by proving (i). By Corollary \ref{res}, we have that $y^K=(Ky)_{\pi}$ and hence
$$\frac{|\{z \in(Ky)_{\pi}|[x,z]=1 \}|}{|(Ky)_{\pi}|}=\frac{|y^{K}\cap C_{G}(x)|}{|y^{K}|}.$$

Now, we have that $y^K\cap C_{G}(x)=(Ky)_{\pi}\cap C_G(x)$ and since $Ky\cap C_G(x)=C_K(x)y$, we deduce that $y^K\cap C_{G}(x)=(C_K(x)y)_{\pi}$. Now, applying again Corollary \ref{res}, we have that $(C_K(x)y)_{\pi}=y^{C_K(x)}$. Since $C_{C_K(x)}(y)=C_K(x)\cap C_K(y)$, we deduce that
$$\frac{|y^{K}\cap C_{G}(x)|}{|y^{K}|}=\frac{|C_{K}(x):C_{K}(x) \cap C_{K}(y)|}{|K:C_{K}(y)|}.$$

 If this proportion is $1$, then $[y,K]\leq C_{G}(x)$ and hence, by the Three Subgroups Lemma, we have that $1=[y,[x,K]]=[y,K]$, which is impossible. Thus, the proportion is strictly smaller than $1$, and hence applying Lemma \ref{inequality} with $P=\langle x \rangle$ and $L=C_{K}(y)$, we have that this proportion is at most $\frac{1}{p+1}$ and the result follows.

 Now, we prove (ii).  For each $i \in \{0,1\ldots p-1\}$ we define $a_{i}=|(Kx^{i}y)_{\pi}\cap C_G(x)|$ and $b_{i}=|(Kx^{i}y)_{\pi}|$. Then
 $$\frac{|\{z \in(Ly)_{\pi}|[x,z]=1 \}|}{|(Ly)_{\pi}|}=\frac{\sum_{i=0}^{p-1}a_{i}}{\sum_{i=0}^{p-1}b_{i}}.$$
 
Let $i \in \{0,1\ldots p-1\}$. If $[x^{i}y,K]\not=1$, then, by part (i), we have that  $a_{i}\leq \frac{b_{i}}{p+1}$. Thus, if $[x^{i}y,K]\not=1$ for all $i \in \{0,1\ldots p-1\}$, then 
$$\frac{\sum_{i=0}^{p-1}a_{i}}{\sum_{i=0}^{p-1}b_{i}} \leq \frac{\frac{1}{p+1}\sum_{i=0}^{p-1}b_{i}}{\sum_{i=0}^{p-1}b_{i}} = \frac{1}{p+1}<\frac{1}{p}$$ and the result follows. 
Therefore, we may assume that there exists $i \in \{0,1\ldots p-1\}$ such that  $[x^{i}y,K]=1$.  Without loss of generality, we may assume that $i=0$, and hence $[y,K]=1$. Thus, $y^{K}=\{y\}$ and hence $a_{0}=b_{0}=1$. Now, since $[y,K]=1$ and $[x,K]=K$, we have that $[x^{i}y,K]\not=1$ for all $i \in \{1,\ldots ,p-1\}$. Thus, we have that  $a_{i}\leq \frac{b_{i}}{p+1}$ for all $i\in \{1,\ldots ,p-1\}$ and hence $$\frac{1+\sum_{i=1}^{p-1}a_{i}}{1+\sum_{i=1}^{p-1}b_{i}}\leq \frac{1+\frac{1}{p+1}\sum_{i=1}^{p-1}b_{i}}{1+\sum_{i=1}^{p-1}b_{i}}=\frac{p}{(p+1)(1+\sum_{i=1}^{p-1}b_{i})}+\frac{1}{p+1}.$$
We also observe that $1+\sum_{i=1}^{p-1}b_{i}=|(Ly)_{\pi}|$. Now, given $i \in \{1,\ldots, p-1\}$ we have that $x^iy\in (Kx^iy)_{\pi}$ and that $x^iy$ commutes with $x$. Thus, $b_i\geq (p+1)a_i\geq p+1$ and hence $|(Ly)_{\pi}|\geq 1+(p-1)(p+1)=p^2$. Therefore, 
$$\frac{1+\sum_{i=1}^{p-1}a_{i}}{1+\sum_{i=1}^{p-1}b_{i}}\leq \frac{p}{(p+1)|(Ly)_{\pi}|}+\frac{1}{p+1}\leq \frac{p}{(p+1)p^2}+\frac{1}{p+1}=\frac{1}{p}$$
and the result follows.
\end{proof}
\end{pro}

\begin{thm}\label{easy}
Let $\pi$ be a set of primes, let $G$ be a finite group and let $x \in G$ be an element of order $p$, where $p \in \pi$. If there exists a normal $\pi'$-subgroup $K$  of $G$ such that $[x,K]\not=1$, then
$$\frac{|C_{G}(x)_{\pi}|}{|G_{\pi}|}\leq \frac{1}{p}.$$
\begin{proof}
Since $KC_{G}(x)$ contains all elements in  $G$ commuting with $x$ we may assume that $G=KC_G(x)$. Since $(o(x),K)=1$, we have that $[x,K]=[x,[x,K]]$ (see \cite[Lemma 4.29]{Isaacs}). Moreover,  $[x,K]$  is normal in $G$ and  $[x,K]$ is a  $\pi'$-subgroup (since $[x,K]\leq K$). Thus, we may replace $K$ by $[x,K]$ and hence, we may assume that $K=[x,K]$. Let $L=\langle K, x \rangle$. It is easy to see that every $\pi$-element of $G$ lies in  $Ly$ for some $y\in G_{\pi}$. Therefore, it suffices to prove that the proportion of $\pi$-elements in $Ly$ commuting with $x$ is at most $\frac{1}{p}$ for every  $y\in G_{\pi}$. That is, we aim to prove that
\begin{equation} \label{eq:1}
\frac{|\{z \in (Ly)_{\pi}| [z,x]=1 \}|}{|(Ly)_{\pi}|}\leq \frac{1}{p}
\end{equation}
for every $y \in G_{\pi}$.

 If no $\pi$-element of $Ly$ commutes with $x$, then left-hand side of inequality (\ref{eq:1}) is $0$ and hence  inequality (\ref{eq:1}) holds trivially. Let us assume that there exists a $\pi$-element in $Ly$ commuting with $x$. Thus, we may assume that $[x,y]=1$ and hence,  inequality (\ref{eq:1}) holds   by Proposition \ref{first} (ii).
\end{proof}
\end{thm}

Finally, we restate and prove Theorem B.

\begin{thm}
Let $\pi$ be a set of primes and let $G$ be  a finite group such that $\Fa(G)$ is a $\pi'$-group. If $x$ is a non-trivial $\pi$-element of $G$, then 
$$\frac{|C_{G}(x)_{\pi}|}{|G_{\pi}|} \leq \frac{1}{q}.$$
where $q$ is the largest prime dividing $o(x)$.
\begin{proof}
Since $C_{G}(x)\leq C_{G}(x^{\frac{o(x)}{q}})$, we may assume that $x$ is an element of order $q$. Now, we know that  $\Fa(G)$ is a normal subgroup such that $C_{G}(\Fa(G))\leq \Fa(G)$ (see \cite[Theorem 9.8]{Isaacs}). By hypothesis, we have that $(o(x),|\Fa(G)|)=1$, and hence $[x,\Fa(G)]\not=1$.  Thus, applying Theorem \ref{easy} with $K=\Fa(G)$, we have that $$\frac{|C_{G}(x)_{\pi}|}{|G_{\pi}|}\leq \frac{1}{q},$$
and the result follows.
\end{proof}
\end{thm}

\section{Proof of Theorem A}\label{TheoremA}

In this section, we prove  Theorem A.  We remark that the proof of Theorem A relies on Theorem \ref{fpr1} and Corollary \ref{fpr2}. Before beginning with the proof, we have to prove some preliminary results, which provide information on the proportion of $\pi$-elements commuting with a fixed automorphism of  order $p \in \pi$ of a central product of  quasisimple groups.

\begin{pro}\label{nuevo37}
Let $K$ be a central product of quasisimple groups with $O_{\pi}(K)=1$ and assume that the simple quotients of all components of $K$ are isomorphic. Let $x,y \in \Aut(K)$ be nontrivial $\pi$-elements such that $x$  acts nontrivially on the set of components of $K$. Then the proportion of elements in $y^K$ commuting with $x$ is at most $\frac{1}{p+1}$, where $p$ is the smallest prime dividing $o(x)$.   
\begin{proof}
If no element of $y^K$ commutes with $x$, then the result holds trivially. Thus, we may assume that there exists an element in $y^K$ commuting with $x$ and hence, replacing $y$ by an appropriate $K$-conjugate of $y$, we may assume that $[x,y]=1$. Let $K_1,\ldots,K_t$ be the components of $K$ and let $L$ be the simple group such that $K_i/Z(K_i)\cong L$ for all $i \in \{1,\ldots , t\}$. Since $x$ does not act trivially on the components, we have that $t\geq p>1$. Set  $J=K/Z(K)\cong L^t$. Since $\Aut(K)\cong \Aut(J)$, then we can view $x$ and $y$ as automorphisms of $J$. Thus, it suffices to prove that the proportion of elements in $y^J$ commuting with $x$  (considered as automorphisms of $J$)  is at most  $\frac{1}{p+1}$.  We define $G=\langle J,x,y\rangle \leq \Aut(J)$. We observe that $J$ is the unique minimal normal subgroup of $G$ and hence $J$ is the socle of $G$. It is not hard to see that  for $z\in J$, $y^{z} \in C_{G}(x)$ if and only if $x^{z^{-1}} \in C_{G}(y)$. Therefore, 
 $$\frac{|y^{J}\cap C_{G}(x)|}{|y^{J}|}=\frac{|x^{J}\cap C_{G}(y)|}{|x^{J}|}.$$
 
 Now, let $H$ be a maximal subgroup of $G$ containing $C_G(y)$. Thus, $G$ acts primitively on the set of cosets $\Omega=G/H$ and hence, we have that
 $$\frac{|x^{J}\cap C_{G}(y)|}{|x^{J}|}\leq  \frac{|x^G\cap H|}{|x^G|} =\fpr(x,\Omega).$$
 Thus, it suffices to prove that $\fpr(x,\Omega)\leq \frac{1}{p+1}$.
 




Assume first that $G\leq L \wr S_t$ with its product action on $\Omega=\Gamma^t$, where $\Gamma$ is a set such that $L\leq \Sym(\Gamma)$.   Thus, $x=nh$ where $n\in \Aut(L)^t$ and $h\in S_t$. Since $x$ does not fix any component, we have that $h\not=1$. Since $\fpr(x,\Omega)\leq \fpr(x^k,\Omega)$ for any $k \geq 1$, we may replace $x$ by any of its powers. Therefore, we may assume that $o(h)=r$ for some prime $r$ dividing $o(x)$ (in particular, $r\geq p$). Since we are assuming that $o(h)=r$, then, reasoning as in the proof of \cite[Theorem 6.1]{FPR}, it is possible to prove that 
 $$\fpr(x,\Omega)\leq \frac{1}{r+1}\leq \frac{1}{p+1}$$
  and the result follows in this case.

 Assume now that the action of $G$ on $\Omega$ is not a product type action. We claim that  $(G,\Omega)$ is none of the exceptions (i), (ii) or (iii) in Theorem \ref{fpr1}. By our assumption, we are not in case (iii) of Theorem \ref{fpr1}. Now, since the socle of $G$  is $J\cong L^t$ with $t>1$, we deduce that we have that $G$ is not almost simple and hence we are not in case (i) of Theorem \ref{fpr1}. Finally, since the socle of $G$ is $J$, which is a non-solvable group, we have that we are not  in case (ii) of Theorem \ref{fpr1}. Thus, the claim follows.

 
 
 Now, we observe that $x^{\frac{o(x)}{p}}$ is an element of order $p$ in $G$. Therefore, since $(G, \Omega)$ is none of the exceptions (i), (ii) or (iii) in Theorem \ref{fpr1}, we have that 
 $$\fpr(x,\Omega)\leq\fpr(x^{\frac{o(x)}{p}},\Omega)\leq \frac{1}{p+1}$$
 and the result follows.
\end{proof}
\end{pro}

Observe that, in the last paragraph, it may occur that $x^{\frac{o(x)}{p}}$ fixes all components of $K$. However, that does not change the argument. Before continuing, we recall and prove some results that we will need.

\begin{pro}[Lemma 3.9 of \cite{Orig}]\label{nuevo39}
Let $G$ be an almost simple group with socle $J$ and assume $J$ is not isomorphic
to an alternating group. Let $p$ be a prime divisor of $|J|$ and suppose $x \in G$ has order $p$. Then there exists an element $y \in J$ of order $p$ such that
$$\frac{|y^G \cap C_G(x)|}{|y^G|}\leq \frac{1}{p+1}.$$
\end{pro}

The next result is a refinement of  \cite[Lemma 3.7]{Orig}.

\begin{lem}\label{nuevo36}
Let $\pi$ be a set of primes, let $G$ be a finite group,  let $x \in G_{\pi}\setminus O_{\pi}(G)$ and let $$D=\{y \in G | \langle y \rangle \text{ is conjugate to } \langle x \rangle \}.$$  
Then $|D|\geq p^2-1$, where $p$ is the smallest prime dividing $o(x)$.
\begin{proof}
If $x$ is a $q$-element for some prime $q \in \pi$, then the result follows from \cite[Lemma 3.7]{Orig}.  Thus, we may assume that $o(x)$ is divisible by at least two different primes. Let us define 
$$T=\{y \in \langle x \rangle|\langle y \rangle=\langle x \rangle\}.$$ 
It is easy to see that $T\subseteq D$ and $|T|=\varphi(o(x))$, where $\varphi$ denotes Euler’s totient function. Thus, if $\varphi(o(x))\geq p^2-1$, then the result holds.

Assume first that $p>2$. If $o(x)$ has at least three prime divisors counting multiplicities, then $\varphi(o(x))\geq (p-1)^3\geq p^2-1$. Thus, we may assume that $o(x)$ is divisible by at most two primes counting multiplicities. Since $x$ is not a $q$-element, we have that $o(x)=pt$, where $p<t$ are two odd primes (in particular, $t\geq p+2$). Therefore, $\varphi(o(x))=(p-1)(t-1)\geq (p-1)(p+1)=p^2-1$.

Assume now that $p=2$. We may assume that $\varphi(o(x))\leq 2$ and hence $o(x)\in \{2,3,4,6\}$. Since $x$ cannot be a $q$-element,  we deduce that $o(x)=6$. In this case, $\{2,3\}\subseteq \pi$ and $T=\{x,x^5\}$. Assume that $x^g\in T$ for every $g \in G$. In this case, $\langle x \rangle $ is normal in $G$ and hence $x \in O_{\pi}(G)$, which is a contradiction. Thus, there exists $g \in G$ such that $x^g \not \in T$. Then $\{x,x^5,x^g\}\subseteq D$ and hence $|D|\geq 3$. 
\end{proof}
\end{lem}

\begin{pro}\label{nuevo310}
Let $K$ be a  quasisimple group with $O_{\pi}(K)=1$. Let $J=K/Z(K)$, and let $x \in \Aut(K)_{\pi}$ such that $x$ has order $p\in \pi$. Assume that $J$ is not an alternating group. If $y \in \Aut(K)$ is a $\pi$-element, then the proportion of elements in $(Ky)_{\pi}$ commuting with $x$ is at most $\frac{1}{p}$. 
\begin{proof}
Reasoning as in Proposition \ref{nuevo37}, we may assume that $[x,y]=1$ and view $y$ as an automorphism of $J$. Set $G=\langle J,x,y \rangle$. We will divide the proof in different claims.

\underline{First claim:}  If $(G,p)\not=(O_{n}^{+}(2),2)$ and $z\in G \setminus \{1\}$, then the proportion of elements in $z^K$ commuting with $x$ is at most $1/p$.

 Assume first that $(J,p)$ is not an exception  in Corollary \ref{fpr2}. Let  $H$ be any maximal core-free subgroup of $G$ containing $C_{G}(z)$. Since we are assuming that $(J,p)$ is not an exception in Corollary \ref{fpr2}, we have that
 $$\frac{|z^{J}\cap C_{G}(x)|}{|z^{J}|}=\frac{|x^{J}\cap C_{G}(z)|}{|x^{J}|}\leq \frac{|x^G \cap H|}{|x^G|} = \fpr(x,G/H) \leq \frac{1}{p}.$$
 Reasoning as in Theorem \ref{nuevo37}, we deduce that  the proportion of elements in $z^K$ commuting with $x$ is at most $\frac{1}{p}$.
 
 Now, we prove the case when $(J,p)$ is one of the exceptions of Corollary \ref{fpr2} (ii).  Let $r$ be a prime and let $z\in G$ be an element of order $r$ and let $(G,p,r)\not=(O_{n}^{+}(2),2,2)$. By looking at  \cite[Table 1]{FPR}, we observe that either $p<r$ or one of the following holds:

\begin{itemize}
\item[(i)] $r=p=2$ and $G\not=O_{n}^{+}(2)$.

\item[(ii)] $r<p=q-1$ and $J=\PSL_2(q)$.

\item[(iii)]$(G,p,r)=(\PSp_6(2),3,2)$.

\end{itemize}

 Assume first that we are in  cases (i), (ii) or  (iii). Since we are excluding the case $(G,p)\not=(O_{n}^{+}(2),2)$, we have that, in any of these cases, there exists some maximal subgroup $H$  of $G$ such that $C_G(z)\leq H$ and $\fpr(x,G/H)\leq 1/p$. Thus, 
  $$\frac{|z^{J}\cap C_{G}(x)|}{|z^{J}|}=\frac{|x^{J}\cap C_{G}(z)|}{|x^{J}|}\leq \frac{|x^G \cap H|}{|x^G|} = \fpr(x,G/H) \leq \frac{1}{p},$$
 and the claim follows in this case. Assume now that $p<r$ and let $T$ be a maximal subgroup of $G$ containing $C_G(x)$. Then $\fpr(z,G/T)\leq 1/r$ and hence 
 $$\frac{|z^{J}\cap C_{G}(x)|}{|z^{J}|}\leq \frac{|z^G \cap T|}{|z^G|} = \fpr(z,G/T) \leq \frac{1}{r}<\frac{1}{p}.$$
 
 That proves that the first claim holds when $(J,p)$ is one of the exceptions and $z$ is an element of prime order. Since $C_G(z)\subseteq C_G(z^m)$ and $\fpr(z,G/H)\leq  \fpr(z^m,G/H)$, for all  $z \in G\setminus \{1\}$ and   for all $m\geq 1$, we have that the first claim holds for arbitrary $z \in G\setminus \{1\}$.
 
 We also observe that if $(G,p)\not=(O_{n}^{+}(2),2)$ and $z\in G\setminus \{1\}$ and it is not a $2$-element, then reasoning as above, the proportion of elements in $z^K$ commuting with $x$ is at most $1/3$.

 \underline{Second claim:} There exists a normal subset $D$ of $K_p\setminus \{1\}$  such that $|D|\geq p^2-1$ and the proportion of elements in $D$  commuting with $x$ is at most $1/(p+1)$.  
 
  By Theorem \ref{nuevo39}, we have that there exists $w$ a non-trivial $p$-element of $G$ such that 
$$\frac{|w^J \cap C_G(x)|}{|w^J|} \leq \frac{|w^G \cap C_G(x)|}{|w^G|}\leq \frac{1}{p+1}.$$
Let us identify $w$ with its image on $K$. Thus, if we embed $w$ in $K$ and we set 
 $$D=\{z \in G | \langle z \rangle \text{ is conjugate to } \langle w \rangle \text{ in } K\},$$
then  we have that $D$ is a normal subset of $p$-elements. Therefore, $D=\bigcup_{i=1}^{r} z_{i}^{K}$, where each $z_j$ is conjugate to $w$ in $K$ and hence
 $$\frac{|z_{i}^{J} \cap C_G(x)|}{|z_{i}^{J}|}=\frac{|w^J \cap C_G(x)|}{|w^J|}\leq \frac{1}{p+1}.$$
 Thus, the proportion of elements in $D$ commuting with $x$ is at most $1/p$ and the claim follows. Moreover, by  Lemma  \ref{nuevo36}, we  know that $|D|\geq p^2-1$.


 Now, we prove the result when $(G,p)\not =(O_{n}^{+}(2),2)$. We have that $(Ky)_{\pi}=\bigcup_{i=1}^{t}y_i^{K}$, where  each  $y_{i}^{K}$ is a $\pi$-element. Assume first that $Ky\not=K$.  By the first claim, we have  that the  proportion of elements in $y_{i}^K$ commuting with $x$ is at most $1/p$ for each $i$ and the result follows in this case. Assume now that $Ky=K$. In this case, we have to count the trivial conjugacy class. Let $D$ be the set defined in the second claim. Since $D\subseteq K_{\pi}\setminus \{1\}$ and $D$ is normal in $G$, the we have that $(Ky)_{\pi}=K_{\pi}=\{1\}\cup D\cup (\bigcup_{i=1}^{s}y_{i}^{K}) $, where  each  $y_{i}^{K}$ is a non-trivial $\pi$-element.  By the first claim, we have that the  proportion of elements in $y_{i}^K$ commuting with $x$ is at most $1/p$ for all $i$. Therefore, it suffices to prove that
   $$\frac{|D\cap C_{K}(x)|+1}{|D|+1}\leq \frac{1}{p}.$$ 
  Now, by the second claim, we have that $|D\cap C_{K}(x)|\leq \frac{|D|}{p+1}$. Thus,
    $$\frac{|D\cap C_{K}(x)|+1}{|D|+1}\leq \frac{1+\frac{1}{p+1}|D|}{|D|+1}=\frac{p+1+|D|}{(p+1)(|D|+1)}=$$ 
    $$=\frac{p}{(p+1)(|D|+1)}+\frac{1}{p+1}\leq \frac{1}{(p+1)p}+\frac{1}{p+1} = \frac{1}{p},$$
 where the last inequality holds because $|D|\geq p^2-1$.

 Thus, it only remains to consider the case when $(G,p)=(O_{n}^{+}(2),2)$ (in particular, $2\in \pi$). Assume first that $Ky$ contains no $2$-element. In this case $(Ky)_{\pi}=\bigcup_{i=1}^{t}y_i^{K} $, where each $y_i$ is a $\pi$-element which is not a $2$-element. Thus, applying the comment after the first claim, we deduce that  the proportion of elements in $y_{i}^{K}$ which commute with $x$ is at most $1/3$ for each $i$ and hence, the result follows in this case. Assume now that $Ky$ contains a $2$-element. Replacing $y$ by a $2$-element in $Ky$, we may assume that $y$ is a $2$-element. In this case, we have that $(Ky)_{\pi}=(Ky)_2\cup (\bigcup_{i=1}^{t}y_{i}^{K}) $, where each $y_i$ is a $\pi$-element which is not a $2$-element. Reasoning as before,  the proportion of elements in $y_{i}^{K}$ which commute with $x$ is at most $1/3$ for each $i$. Moreover, by  \cite[Proposition 3.10 (ii) b)]{Orig}, we have that the proportion of elements in $(Ky)_2$ commuting with $x$ is at most $1/2$. Thus, the result follows in this case.
\end{proof}
\end{pro}

Finally, we prove two results that  will allow us to study the case when the group has alternating composition factors.

\begin{pro}\label{SymNonTrans}
Let $p$ be a prime, let $\pi$ be a set of  primes containing $p$ and let $n \geq 5$. If  $x \in S_n$ is an element of order $p$ and $x$ is not a transposition, then 
$$|(A_n)_{\pi}|\geq p |C_{A_n}(x)_{\pi}|,$$
and
$$|(S_n \setminus A_n)_{\pi}|\geq p |(C_{S_n}(x)\setminus A_n)_{\pi}|.$$
Equivalently, the proportion of $\pi$-elements in  $A_n$ and in $S_n\setminus A_n$ commuting with $x$ (the last proportion is defined as $0$ if $2 \not \in \pi$)  is at most $\frac{1}{p}$.
\begin{proof}
 Let $m$ be the size of the support of $x$. We notice that $p$ divides $m$ and that $C_{S_n}(x)=(C_p\wr S_{m/p})\times S_{n-m} \leq (S_p\wr S_{m/p})\times S_{n-m}$. 
 
 Assume first that $m=p$ and that $p\geq5$. In this case, $C_{A_n}(x)=\langle x \rangle \times S_{n-m}$ and hence, $|(C_{A_n}(x))_{\pi}|=p|(A_{n-p})_{\pi}|$ and $|(C_{S_n}(x)\setminus A_n)_{\pi}|=p|(S_{n-p}\setminus A_{n-p})_{\pi}|$. On the other hand, we also have that $|(A_n)_{\pi}|\geq |(A_p)_{\pi}||(A_{n-p})_{\pi}|$ and $|(S_n \setminus A_n)_{\pi}|\geq |(A_p)_{\pi}||(S_{n-p}\setminus A_{n-p})_{\pi}|$. Since $A_p$ does not possess a a normal Sylow $p$-subgroup, we have that $|(A_p)_{p}|\geq p^2$ and hence $|(A_p)_{\pi}|\geq p^2$. Thus,  both inequalities of the proposition hold.


  Assume now that $m>p$ and that $m\geq 5$. Since $A_m$ does not possess a normal Sylow $p$-subgroup, we have that $|(A_m)_{\pi}|\geq p^2$.  In this case, we have the following inequalities
 $$|C_{A_n}(x)_{\pi}|\leq |((S_p\wr S_{m/p}) \cap A_m)_{\pi}||(A_{n-m})_{\pi}|+|((S_p\wr S_{m/p}) \cap (S_m \setminus A_m))_{\pi}||(S_{n-m} \setminus A_{n-m})_{\pi}|$$
  and 
   $$|(C_{S_n}(x)\setminus A_n)_{\pi}|\leq  |((S_p\wr S_{m/p})  \setminus A_m)_{\pi}||(A_{n-m})_{\pi}|+|((S_p\wr S_{m/p}) \cap A_m)_{\pi}||(S_{n-m} \setminus A_{n-m})_{\pi}|.$$
   
Now, we claim that \begin{equation} \label{eq:2}
|((S_p\wr S_{m/p}) \cap A_m)_{\pi}|\leq \frac{|(A_m)_{\pi}|}{p}
\end{equation}
and
 \begin{equation} \label{eq:3}
 |((S_p\wr S_{m/p}) \setminus A_m))_{\pi}|\leq \frac{|(S_m \setminus A_m)_{\pi}|}{p}.
\end{equation}
 
 Let us prove inequality (\ref{eq:2}) (the proof of inequality (\ref{eq:3}) is analogous). Let $(A_m)_{\pi}=\{1\}\cup  (\bigcup_{i=1}^{t}y_{i}^{A_m})$, where each $y_i$ is a non-trivial $\pi$-element of $A_m$. Applying Theorem \ref{fpr1}, we have that 
 $$\frac{|y_{i}^{A_m}\cap (S_p \wr S_{m/p})|}{|y_{i}^{A_m}|}\leq \frac{1}{p+1}$$
 for all $i \in \{1, \ldots, t\}$. Thus, we deduce that 
 $$\frac{|((S_p\wr S_{m/p}) \cap A_m)_{\pi}|}{|(A_m)_{\pi}|}=\frac{1+\sum_{i=1}^{t}|y_{i}^{A_m}\cap (S_p \wr S_{m/p})|}{1+\sum_{i=1}^{t}|y_{i}^{A_m}|}\leq$$
 $$\leq \frac{p+1+\sum_{i=1}^{t}|y_{i}^{A_m}|}{(p+1)(1+\sum_{i=1}^{t}|y_{i}^{A_m}|)}=\frac{p}{(p+1)|(A_m)_{\pi}|}+\frac{1}{p+1}\leq \frac{1}{p},$$
 where the last inequality holds since $|(A_m)_{\pi}|\geq p^2$. Thus, inequality (\ref{eq:2}) holds.

 Now, applying  inequalities (\ref{eq:2}) and (\ref{eq:3}), we have that

  $$|C_{A_n}(x)_{\pi}|\leq \frac{1}{p}(|(A_m)_{\pi}||(A_{n-m})_{\pi}|+|(S_m \setminus A_m)_{\pi}||(S_{n-m} \setminus A_{n-m})_{\pi}|)=$$ $$= \frac{1}{p}|(S_m \times S_{n-m})_{\pi}\cap A_n|\leq \frac{|(A_n)_{\pi}|}{p}.$$
 and 
   $$|(C_{S_n}(x)\setminus A_n)_{\pi}|\leq  \frac{1}{p}(|(S_m \setminus A_m)_{\pi}||(A_{n-m})_{\pi}|+|(A_m)_{\pi}||(S_{n-m} \setminus A_{n-m})_{\pi}|)=$$
$$= \frac{1}{p}  |(S_m\times S_{n-m})_{\pi}\setminus A_n|  \leq \frac{|(S_n\setminus A_n)_{\pi}|}{p}$$
and hence both inequalities of the proposition hold.

Now, it only remains to study the case when $m \in \{3,4\}$.
 
 
 Suppose first that $m=3$. Then $3 \in \pi$ and  $x$ is a $3$-cycle. Without loss of generality, we may assume that $x=(1,2,3)$. In this case, we have that $C_{S_{n}}(x)_{\pi}=\langle (1,2,3) \rangle \times (S_{n-3})_{\pi}$  and hence $|C_{A_{n}}(x)_{\pi}|=3| (A_{n-3})_{\pi}|$ and $|(C_{S_{n}}(x)\setminus A_n)_{\pi}|=3| (S_{n-3}\setminus A_{n-3})_{\pi}| $. Fix $i\in \{4,\ldots n\}$. Let $L(i)$ (respectively, $T(i)$) be the set of elements of the form $(1,2,i)b$ or $(1,i,2)b$ for $b\in (A_n)_{\pi}$ (respectively, $b\in (S_n \setminus A_n)_{\pi}$) fixing $1$,$2$ and $i$. Thus, $L(i)$ are disjoint subsets of elements of $ (A_n)_{\pi}$  of size $2| (A_{n-3})_{\pi}|$,  whose elements do not commute with $x$. Adding the number of $\pi$-elements in $A_{n}$ which commute with $x$, we have that $|(A_n)_{\pi}|\geq(2n-3) | (A_{n-3})_{\pi}|$. Replacing the $L(i)$ by $T(i)$ and reasoning analogously,  we deduce that  $|(S_n\setminus A_n)_{\pi}|\geq(2n-3) | (S_{n-3}\setminus A_{n-3})_{\pi}|$. Therefore,
$$\frac{|C_{A_{n}}(x)_{\pi}|}{|(A_{n})_{\pi}|}\leq  \frac{3}{2n-3},$$
and we observe that $\frac{3}{2n-3}$ is smaller than $1/3$ for $n \geq 6$. For $n=5$ we only have to compute this proportion for  $\pi\in \{\{3\},\{2,3\},\{3,5\}, \{2,3,5\}\}$. Since $C_{A_{5}}(x)= \langle x \rangle$, we have that $|C_{A_{5}}(x)_{\pi}|=3$ for any of the possible $\pi$. Moreover, $|(A_5)_{\pi}|\geq |(A_5)_{3}|=21$. Thus, $\frac{|C_{A_{5}}(x)_{\pi}|}{|(A_5)_{\pi}|}\leq \frac{3}{21}<\frac{1}{3}$ and the result follows in this case. Assume now that $2\in \pi$ (otherwise, $(C_{S_{n}}(x)\setminus A_n)_{\pi}=\varnothing$ and there is nothing to prove). Reasoning as above, we have that 
$$\frac{|(C_{S_{n}}(x)\setminus A_n)_{\pi}|}{|(S_{n}\setminus A_{n})_{\pi}|}\leq  \frac{3}{2n-3},$$
which is smaller than $1/3$ for $n \geq 6$.  For $n=5$ we only have to compute the possibilities $\pi\in \{\{2,3\},\{2,3,5\}\}$. Now, we have that $|(C_{S_5}(x)\setminus A_5)_{\{2,3,5\}}|=|(C_{S_5}(x)\setminus A_5)_{\{2,3\}}|=3$ and hence $\frac{|(C_{S_{5}}(x)\setminus A_5)_{\pi}|}{|(S_5 \setminus A_5)_{\pi}|}\leq \frac{3}{30}<\frac{1}{3}$. Thus, the result holds when $m=3$.

Finally, assume that $m=4$. In this case we may assume that $x=(1,2)(3,4)$ and hence $C_{S_n}(x)=C_{S_4}(x)\times S_{n-4}$ and hence
$$|C_{A_n}(x)_{\pi}|=|(C_{S_n}(x)\setminus A_n)_{\pi}|=4|(S_{n-4})_{\pi}|.$$

 Let $5\leq i<j\leq n $. Let $Z(i,j)$ be the set of elements of $S_n$ of the form $uv$, where $u$ is a double transposition  on $\{1,2,i,j\}$ different from $(1,2)(i,j)$ and $v$ is a $\pi$-element fixing each element of $\{1,2,i,j\}$.  Analogously, let $W(i,j)$ be a set defined as $Z(i,j)$, but replacing $u$ by a $4$-cycle on $\{1,2,i,j\}$ different from $(1,2,i,j)$.  Thus, $W(i,j)$ and $Z(i,j)$ are sets of $\pi$-elements not commuting with $x$ such that $|W(i,j)|=2|(S_{n-4})_{\pi}|$ and $|Z(i,j)|=2|(S_{n-4})_{\pi}|$. Thus, for each choice of $\{i,j\}$, we obtain at least $2|(S_{n-4})_{\pi}|$ different $\pi$-elements in $A_n$ or in $S_n\setminus A_n$ not commuting with $x$. Thus,
 $$|(A_n)_{\pi}|,|(S_n \setminus A_n)_{\pi}|\geq (4+(n-4)(n-5))|(S_{n-4})_{\pi}|.$$
 
Thus, if $n \geq 7$ we have that $\frac{4}{4+(n-4)(n-5)}\leq \frac{1}{2}$ and then, both inequalities of the proposition hold. Therefore, we may assume that $n\in \{5,6\}$. Let $\pi$ be a set of primes  such that $\{2\} \subseteq \pi \subseteq \{2,3,5\}$ and let $n \in \{5,6\}$. Then $C_{S_n}(x)_{\pi}=C_{S_n}(x)_{2}$ and hence
$$\frac{|C_{A_n}(x)_{\pi}|}{|(A_n)_{\pi}|}\leq \frac{|C_{A_n}(x)_{2}|}{|(A_n)_{2}|}\leq \frac{1}{2}$$
and
$$\frac{|(C_{S_n}(x)\setminus A_n)_{\pi}|}{|(S_n\setminus A_n)_{\pi}|}\leq \frac{|(C_{S_n}(x)\setminus A_n)_{2}|}{|(S_n\setminus A_n)_{2}|}\leq \frac{1}{2}.$$ Thus, the result holds when $m=4$.
\end{proof}
\end{pro}

\begin{pro}\label{SymTrans}
Let $n\geq 5$ and let $\pi$ be a set of primes containing $2$. If $x\in S_n$ is a transposition, then 
$$\frac{|C_{S_n}(x)_{\pi}|}{|(S_n)_{\pi}|}<\frac{1}{2}.$$
\begin{proof}
We may assume that $x=(1,2)$. We know that $C_{S_n}(x)=\langle (1,2)\rangle \times S_{n-2}$ and hence $|C_{S_n}(x)_{\pi}|=2|(S_{n-2})_{\pi}|$. Given $j \in \{3,\ldots , n\}$ we define $Z_j=\{\tau \in (S_n)_{\pi}|\tau(1)=j, \tau(j)=1\}$. Now,   we have that $Z_j$ is a set of $\pi$-elements not commuting with $x$. In addition, we have that $Z_i\cap Z_j=\varnothing$ if $i\not=j$. Thus, counting the $\pi$-elements centralizing $x$, we have that $|(S_n)_{\pi}|\geq n|(S_{n-2})_{\pi}|$. Therefore, 
$$\frac{|C_{S_n}(x)_{\pi}|}{|(S_n)_{\pi}|}\leq \frac{2|(S_{n-2})_{\pi}|}{n|(S_{n-2})_{\pi}|}=\frac{2}{n}\leq \frac{2}{5}<\frac{1}{2}$$
and the result follows.
\end{proof}
\end{pro}

Now, we restate and prove Theorem A.

\begin{thm}
Let $G$ be  a finite group and let $\pi$ be a set of primes. If $x \in G_{\pi}\setminus O_{\pi}(G)$, then 
$$\frac{|C_{G}(x)_{\pi}|}{|G_{\pi}|} \leq \frac{1}{p},$$
where $p$ is the smallest prime in $\pi.$
\begin{proof}
First, we may assume that $o(x)=q$, where $q$ is a prime in $\pi$ (in particular, $p\leq q$). We also have that $x \not \in \Fi(G)$. By Lemma \ref{gruponorm}, we may assume that $O_{\pi}(G)=1$.

If $x$ does not centralize $O_{\pi'}(G)$, then the result follows by Theorem \ref{easy}.  Therefore, we may assume that $x$ centralizes $O_{\pi'}(G)$, and hence $x \in C_{G}(\Fi(G))$. It follows that $G$ is a non-solvable group, and that $x$ acts faithfully on $\Fa(G)$. As a consequence, $x$ acts faithfully on $K$, where $K$ is a central product of quasisimple components of $G$. We also observe that every $\pi$-element in $G$ which  commutes with $x$ must be a $\pi$-element in $KC_{G}(x)$ which commutes with $x$. Therefore, we may replace $G$ by the subgroup $KC_{G}(x)$ and assume that $C_G(x)$ acts transitively on the components of $K$.  In particular, the simple quotients of all components of $K$ are isomorphic. Thus, we have two possibilities.

\begin{itemize}
\item [(a)] $x$ acts nontrivially on the set of components of $K$. 

\item [(b)] $x$  normalizes each component of $K$.
\end{itemize}

 Now, we observe that the result will follow once we prove that 
$$\frac{|(Ky)_{\pi}\cap C_{G}(x)|}{|(Ky)_{\pi}|}\leq  \frac{1}{p}$$
 for every non-trivial $y \in (C_{G}(x))_{\pi}$. We may also assume that $G=\langle K,x,y \rangle$ and that  $\langle x,y \rangle$ acts transitively on the components of $K$.

Suppose first that we are in case a). Let $z \in (Ky)_{\pi}$ with $z\not=1$. By Proposition \ref{nuevo37}, we have that $$\frac{|z^{K}\cap C_{G}(x)|}{|z^{K}|} \leq \frac{1}{q+1}\leq \frac{1}{p+1}.$$ 
Thus, if $Ky\not=K$, then, expressing $(Ky)_{\pi}$ as a union of conjugacy classes, we have that the proportion of $\pi$-elements in  $Ky$ which commute with $x$ is at most $\frac{1}{p+1}$. Therefore, we may assume that $y \in K$ and hence, if we set  $D=\{z \in K| \langle z \rangle \text{ is conjugate to } \langle y \rangle\}$, then by Lemma \ref{nuevo36}, we have that $|D|\geq p^{2}-1$. Since $D \subseteq K_{\pi}$, we have that $|K_{\pi}|\geq p^{2}$. Now, expressing $K_{\pi}$ as union of $K$-conjugacy classes and reasoning as in the proof of Proposition \ref{first}, we deduce that the proportion of $\pi$-elements in $K$ commuting with $x$ is at most $1/p$.

 Suppose now that we are in case b). It follows that $y$ must act  transitively on the components of $K$. If $K$ has more than two  components, then we can reduce the problem to the previous case (exchanging  $x$ by $y$). Thus, we may assume that $K$ is quasisimple. By Propositions \ref{nuevo310} and \ref{SymNonTrans}, the result holds unless when $K/Z(K)=A_n$ and  that $x$ acts as a transposition. Assume now that $K/Z(K)=A_n$ and  that $x$ acts as a transposition. Thus, if we set $L=\langle K,x \rangle$, then,  by Proposition \ref{SymTrans}, we deduce that the proportion of $\pi$-elements in the coset $Ly$ commuting with $x$ is at most $1/2$. Thus, the proportion of $\pi$-elements in the coset $Ky$ commuting with $x$ is at most $1/2$.
\end{proof}
\end{thm}

As we mentioned in the Introduction, it remains as an open question whether Theorem A holds with $\leq 1/q$, where $q$ is the largest prime dividing $o(x)$. 

\section{Results on $\Pro_{\pi}(G)$}\label{TheoremC}

In this section, we will use Theorem A to prove Theorem C.  If $G$ possesses a normal and abelian Hall $\pi$-subgroup, then every pair of $\pi$-elements of $G$ commute and hence $\Pro_{\pi}(G)=1>\frac{p^{2}+p-1}{p^{3}}$. Thus, we will only prove the converse assertion of Theorem C. We recall a pair of results that we will need.

\begin{pro}[Theorem 5.2 of \cite{Orig}]\label{lem1}
Let $G$ be a non-abelian finite group, and let $p$ be the smallest prime dividing $|G|$. Then $\Pro(G) \leq  \frac{p^{2}+p-1}{p^{3}}$. 
\end{pro}

\begin{lem}\label{PIelem}
Let $G$ be a finite group and let $\pi$ be a set of primes such that $G$ does not possess a normal Hall $\pi$-subgroup. Then $|G_{\pi}|\geq p^2$, where $p$ is the smallest prime in $\pi$.
\begin{proof}
Since $G$ does not possess a normal Hall $\pi$-subgroup, we have that there exists $q \in \pi$ such that $G$ does not possess a normal Sylow $q$-subgroup. 

Let $Q$ be a Sylow $q$-subgroup of $G$. If $|Q|\geq q^2$, then $|G_{\pi}|\geq |Q|\geq q^2\geq p^2$ and the result follows. Assume now that $|Q|=q$. Then all Sylow $q$-subgroups intersect trivially and, by Sylow Theorems, we have that $G$ possesses at least $q+1$ Sylow $q$-subgroups. Thus, $|G_{\pi}|\geq |G_q|\geq q^2\geq p^2$ and the result follows.
\end{proof}
\end{lem}

Now, we proceed to restate and prove Theorem C.

\begin{thm}
Let $G$ be a finite group, let $\pi$ be a set of primes and let $\pi$ be the smallest prime in $\pi$. If $\Pro_{\pi}(G)> \frac{p^{2}+p-1}{p^{3}}$, then $G$ has a normal and abelian Hall  $\pi$-subgroup.
\begin{proof}
First, we  prove that $G$ has a normal Hall $\pi$-subgroup by induction on $|G|$.

Suppose first that  $O_{\pi}(G)>1$. By Lemma \ref{gruponorm}, $\Pro_{\pi}(G/O_{\pi}(G)) \geq \Pro_{\pi}(G)$, and $|G/O_{\pi}(G)|<|G|$. Thus, by induction, $G/O_{\pi}(G)$ has a normal Hall $\pi$-subgroup and hence $G$ has a normal Hall $\pi$-subgroup.

Thus, we may assume that $O_{\pi}(G)=1$. Then,  
$$|\{(x,y)|x,y\in G_{\pi}, xy=yx\}|=|\{(x,y)|x\in G_{\pi},y \in C_{G}(x)_{\pi}\}|= $$
$$=|G_{\pi}|+\sum_{x \in G_{\pi}\setminus \{1\}}|C_{G}(x)_{\pi}|.$$

For every $x \in G_{\pi}\setminus \{1\}$ we have that $x \not \in O_{\pi}(G)$. Thus, applying Theorem A, we have that
$$\frac{|C_{G}(x)_{\pi}|}{|G_{\pi}|}\leq \frac{1}{p},$$
and hence
$$\Pro_{\pi}(G)=\frac{1}{|G_{\pi}|}(1+\sum_{x \in G_{\pi}\setminus \{1\}}\frac{|C_{G}(x)_{\pi}|}{|G_{\pi}|}) \leq \frac{1}{|G_{\pi}|}(1+ \frac{|G_{\pi}|-1}{p}).$$

So, if we consider the function
$$f_{p}(x)=\frac{1}{x}(1+(x-1)\frac{1}{p})=\frac{1}{x}(1-\frac{1}{p})+\frac{1}{p},$$
then we have that $\Pro_{\pi}(G) \leq f_{p}(|G_{\pi}|)$. We observe  that  $f_p(x)$ is a decreasing function.

Assume now that   $G$ does not have a normal Hall $\pi$-subgroup. By Lemma \ref{PIelem},  we have that  $|G_{\pi}|\geq p^{2}$ and hence
$$\Pro_{\pi}(G) \leq  f_{p}(|G_{\pi}|) \leq f_{p}(p^{2})=\frac{1}{p^{2}}(1-\frac{1}{p})+\frac{1}{p}= \frac{p^{2}+p-1}{p^{3}},$$
which is a contradiction.

Thus, $G$ possesses a normal Hall $\pi$-subgroup, or equivalently, $O_{\pi}(G)$ is a Hall $\pi$-subgroup of $G$ and $O_{\pi}(G)=G_{\pi}$. Now, we have that
$$\Pro(O_{\pi}(G))=\frac{|\{(x,y) \in G_{\pi}\times G_{\pi}|xy=yx\}|}{|G_{\pi}|^{2}}=\Pro_{\pi}(G) >\frac{p^{2}+p-1}{p^{3}}.$$
Therefore, $O_{\pi}(G)$ is abelian by  Proposition \ref{lem1}.
\end{proof}
\end{thm}



\begin{thebibliography}{99}



\bibitem[BG22]{FPR} T. C. Burness, R. M. Guralnick, \rm{Fixed point ratios for finite primitive groups and applications},  \textit{Adv. Math.} \textbf{411} (2022).

\bibitem[BGMN23]{Orig} T. C. Burness, R. M. Guralnick, A. Moretó, G. Navarro, \rm{The probability that two $p$-elements commute}, \textit{Alg. Number Th. } \textbf{17} (2023), 1209-1229.




\bibitem[GR06]{GR} R.\,M. Guralnick and G.\,R. Robinson, \rm{On the commuting probability in finite groups}, \textit{J. Algebra} \textbf{300} (2006), 509-528.

\bibitem[Gus73]{Gustafson} W. H. Gustafson, \rm{What is the probability that two elements commute?}, \textit{Amer. Math. Monthly} \textbf{80} (1973), 1031-1034.


\bibitem[HMM23]{HMM}  N. N. Hung,  A. Maróti,   J. Martínez, \rm{Conjugacy classes of $\pi$-elements and nilpotent/abelian Hall $\pi$-subgroups},  \emph{Pacific J. Math.} \textbf{323} (2023), 185-204.

\bibitem[Isa08]{Isaacs} I. M. Isaacs, \textit{Finite Group Theory}, Amer. Math. Soc. Providence, Rhode Island, 2008.

\bibitem[Jo69]{J} K. S. Joseph, \textit{Commutativity in non-abelian groups}. Ph.D. thesis, UCLA, 1969

\bibitem[MN14]{MN} A. Maróti, H. N. Nguyen, \rm{On the number of conjugacy classes of $\pi$-elements in finite groups}, \textit{Arch. Math.} \textbf{102} (2014), 101-108.

\bibitem[Ne89]{N} P.\,M. Neumann, Two combinatorial problems in group theory, \textit{Bull. London Math. Soc.} \textbf{21} (1989), 456-458

\end{thebibliography}
\end{document}